\documentclass[11pt]{article}
\input amssym.def
\input amssym.tex

\def\slashedfrac#1#2{\hbox{\kern.1em %
 \raise.5ex\hbox{\the\scriptfont0 #1}\kern-.11em %
 /\kern-.15em\lower.25ex\hbox{\the\scriptfont0 #2}}}

\def\@regmark{{\ooalign
{\hfil\raise.07ex\hbox{\tiny R}\hfil\crcr{$\scriptstyle\bigcirc$}}}}
\def\regmark{\raise1ex\hbox{\@regmark}}

\newcommand{\tht}{\theta}
\newcommand{\ra}{\rightarrow}

\newcommand{\CC}{{\Bbb C}}

\newcommand{\FF}{{\Bbb F}}

\newcommand{\RR}{{\Bbb R}}

\newcommand{\ZZ}{{\Bbb Z}}

\newcommand{\sF}{{\cal F}}

\newcommand{\sL}{{\cal L}}

\newcommand{\sO}{{\cal O}}
\newcommand{\sP}{{\cal P}}

\newcommand{\Om}{\Omega}

\newcommand{\df}{\displaystyle\frac}

\def\slashedfrac#1#2{\hbox{\kern.1em %
 \raise.5ex\hbox{\the\scriptfont0 #1}\kern-.11em %
 /\kern-.15em\lower.25ex\hbox{\the\scriptfont0 #2}}}

\newcommand{\btd}{\blacktriangledown}
\newcommand{\hsp}{\hspace*{\parindent}}

\newcommand{\eeq}{\end{equation}}
\newcommand{\beql}[1]{\begin{equation}\label{#1}}

\newcommand{\Blb}[1]{\begin{array}{l}
\lceil \\[-.04in] \mbox{\,\rule{.005in}{#1in}} \\[-.13in] \lfloor
\end{array} \hspace*{-.05in}}
\newcommand{\Brb}[1]{\hspace*{-.1in}\begin{array}{c}
\rceil \\[-.04in] \mbox{\hspace*{.005in}\rule{.005in}{#1in}} \\[-.13in] \rfloor
\end{array}}

\makeatletter
\def\section{\@startsection {section}{1}{\z@}{-3.5ex plus -1ex minus 
 -.2ex}{2.3ex plus .2ex}{\normalsize\bf}}
\makeatother
\makeatletter
\def\subsection{\@startsection {subsection}{1}{\z@}{-3.5ex plus -1ex minus
 -.2ex}{2.3ex plus .2ex}{\normalsize\bf}}

\def\@sect#1#2#3#4#5#6[#7]#8{\ifnum #2>\c@secnumdepth
     \def\@svsec{}\else 
     \refstepcounter{#1}\edef\@svsec{\csname the#1\endcsname.\hskip .75em }\fi
     \@tempskipa #5\relax
      \ifdim \@tempskipa>\z@ 
        \begingroup #6\relax
          \@hangfrom{\hskip #3\relax\@svsec}{\interlinepenalty \@M #8\par}%
        \endgroup
       \csname #1mark\endcsname{#7}\addcontentsline
         {toc}{#1}{\ifnum #2>\c@secnumdepth \else
                      \protect\numberline{\csname the#1\endcsname}\fi
                    #7}\else
        \def\@svsechd{#6\hskip #3\@svsec #8\csname #1mark\endcsname
                      {#7}\addcontentsline
                           {toc}{#1}{\ifnum #2>\c@secnumdepth \else
                             \protect\numberline{\csname the#1\endcsname}\fi
                       #7}}\fi
     \@xsect{#5}}
\makeatother

\begin{document}
\begin{center}
{\Large {\bf A Group-Theoretic Framework for the Construction \\
\vspace{.1\baselineskip}
of Packings in Grassmannian Spaces }} \\
\vspace{\baselineskip}
{\em A. R. Calderbank, R. H. Hardin, E. M. Rains,} \smallskip \\
{\em P. W. Shor and N. J. A. Sloane} \\
\vspace{.5\baselineskip}
Information Sciences Research \\
AT\&T Labs - Research \\
180 Park Avenue, Room C233 \\
Florham Park, NJ 07932-0971 \\
\vspace{\baselineskip}
\vspace{1.5\baselineskip}
{\bf ABSTRACT}
\vspace{.5\baselineskip}
\end{center}
\setlength{\baselineskip}{1.5\baselineskip}
\hspace*{.25in}
By using totally isotropic subspaces in an orthogonal space $\Om^+ (2i,2)$,
several infinite families of packings of $2^k$-dimensional subspaces
of real $2^i$-dimensional space
are constructed, some of which are shown to be optimal packings.
A certain Clifford group underlies the
construction and links this problem with Barnes-Wall lattices,
Kerdock sets and quantum-error-correcting codes.
\clearpage
\large\normalsize
\renewcommand{\baselinestretch}{1}
\thispagestyle{empty}
\setcounter{page}{1}
\begin{center}
{\Large {\bf A Group-Theoretic Framework for the Construction \\
\vspace{.1\baselineskip}
of Packings in Grassmannian Spaces }} \\
\vspace{\baselineskip}
{\em A. R. Calderbank, R. H. Hardin, E. M. Rains,} \smallskip \\
{\em P. W. Shor and N. J. A. Sloane} \\
\vspace{.5\baselineskip}
Information Sciences Research \\
AT\&T Labs - Research \\
600 Mountain Avenue \\
Murray Hill, NJ 07974 \\
\vspace{1.5\baselineskip}
\vspace{1.5\baselineskip}
\end{center}
\setlength{\baselineskip}{1.5\baselineskip}
\section{Introduction}
\label{sec1}
\hsp
The central problem is to arrange $N$ $n$-dimensional subspaces of $\RR^m$ so they
are as far apart as possible.
Numerous constructions and bounds were given in
\cite{grass1}, \cite{grass2} (see also \cite{grass0}).
In the present paper we give an algebraic framework for
constructing such arrangements that explains all the examples
constructed or conjectured in \cite{grass2}.

The two main constructions obtained by these methods are stated in Theorems~1 and 2. Theorem~3 describes an unrelated construction which yields another infinite family of optimal packings. Table~I
at the end of the paper summarizes the parameters of the packings obtained in
dimensions up to 128.

$G(m,n)$ will denote the Grassmannian space of $n$-dimensional
subspaces of $\RR^m$.
We shall refer to the elements of $G(m,n)$
as $n$-dimensional
planes, or simply {\em planes}.
For reasons discussed in \cite{grass1}, we define the {\em distance}
between two $n$-dimensional planes
$P,Q$ in $\RR^m$ by
\beql{eq1}
d(P,Q) = \sqrt{\sin^2  \tht_1 + \ldots + \sin^2 \tht_n} ~,
\eeq
where $\tht_1 , \ldots , \tht_n$ are the principal
angles$^1$\footnotetext[1]{\cite{GVL89}, p.~584.}
between $P$ and $Q$.
For given values of $m,n,N$ we wish to find the best
packings of $N$ planes in $G(m,n)$, that is,
subsets
$\sP = \{ P_1 , \ldots , P_N \} \subset G(m,n)$ such that
$d( \sP ) = \min_{i \neq j} d(P_i , P_j )$
is maximized ($d(\sP )$ is called the {\em minimal distance}
of the packing).
We refer to \cite{grass1} for applications and earlier references.

It was shown in \cite{grass1} that $G(m,n)$ equipped with the metric
(\ref{eq1}) has an
isometric embedding in $\RR^D ,~ D = (m-1) (m+2)/2$,
obtained by representing each plane $P \in G(m,n)$ by the
orthogonal projection from $\RR^m$ to $P$.
If $A$ is an $n \times m$ generator matrix for $P$, whose rows
are orthonormal vectors spanning $P$, then the projection is represented by the
matrix $\Pi_P = A^t A$, where $t$ denotes transposition.
$\Pi_P$ is an $m \times m$ symmetric idempotent matrix with
trace $n$, and so lies in $\RR^D$.
All such $\Pi_P$ for $P \in G (m,n)$ lie on a sphere of radius
$\sqrt{n(m-n)/m} $ in $\RR^D$.
Furthermore, if two planes $P,Q \in G(m,n)$ are represented by
projection matrices $\Pi_P , \Pi_Q$ then
\beql{eq2}
d^2 (P,Q) = n - \mbox{trace}~~ \Pi_P~\Pi_Q = \df{1}{2} ||
\Pi_P - \Pi_Q ||^2 ~,
\eeq
where$^2$\footnotetext[2]{\cite{GVL89}, p.~56} $||~~||$ denotes the $L_2$ or
Frobenius norm of a matrix (\cite{grass1}, Theorem~5.1).

It follows from this embedding that if $\sP$ is a packing of $N$ planes in
$G(m,n)$ then
\beql{eq3}
d(\sP)^2  \leq \df{n(m-n)}{m} ~~\df{N}{N-1}
\eeq
(the ``simplex bound'').
Equality requires $N \leq D+1 = {{m+1}\choose{2}}$,
and occurs if and only if the $N$ points in $\RR^D$ corresponding to
the planes form a regular equatorial simplex
(\cite{grass1}, Cor.~5.2).
Also, for $N > D+1$,
\beql{eq4}
d( \sP)^2 \leq \df{n(m-n)}{m}
\eeq
(the ``orthoplex bound'').
Equality requires $N \leq 2D = (m-1)(m+2)$, and occurs if the
$N$ points form a subset of the $2D$ vertices of a regular
orthoplex (generalized octahedron).
If $N = 2D$ this condition is also necessary
(\cite{grass1}, Cor.~5.3).
\section{The algebraic framework}
\label{sec2}
\hsp
The following machinery was used in \cite{CCKS}
to construct Kerdock sets, among other things, and in \cite{QC1},
\cite{QC2} to construct quantum error-correcting codes.
The starting point is the standard method of associating a finite
orthogonal space to an extraspecial 2-group, as described for
example in \cite{Asch}, Theorem~23.10,
or \cite{Huppert}, Theorem~13.8.

The end result will be the construction of various packings of
$n$-spaces in
a parent space $V = \RR^m$, where $m = 2^i$.
As basic vectors for $V$ we use $e_u$, $u \in U = \FF_2^{^i}$.
The constructions will involve certain subgroups of the
real orthogonal group $\sO = O(V, \RR )$.

For $a,b \in U$ we define transformations
$X(a) \in \sO$, $Y (b) \in \sO$ by
$$
X(a): e_u \ra e_{u+a} ~, \quad
Y(b): e_u \ra (-1)^{b \cdot u} e_u ~, ~
u \in U~,
$$
where the dot indicates the usual inner product in $U$.
Then $X = \langle X(a) : a \in U \rangle$,
$Y = \langle Y(b): b \in U \rangle$ are elementary
abelian subgroups of $\sO$ of order $2^i$, and
$E = \langle X,Y \rangle \subset \sO$ is an
extraspecial 2-group$^3$\footnotetext[3]{\cite{Huppert}, p.~349} of order
$2^{2i+1}$
(\cite{CCKS}, Lemma~2.1).
The elements of $E$ have the form $\pm X(a) Y (b)$,
$a,b \in U$, and satisfy
$$
Y(b) X(a) = (-1)^{a \cdot b}
X(a)Y(b)~,
$$

$$
(-1)^s X(a) Y(b) (-1)^{s'}
X(a') Y(b') =
(-1)^{a' \cdot b + s+s'}
X(a+a') Y(b+b') ~.
$$

The center $\Xi (E)$ of $E$
is $\{ \pm I\}$, and $\bar{E} = E/ \Xi (E)$ is an elementary
abelian group of order $2^{2i}$ whose elements
can be denoted by
$\bar{X}(a) \bar{Y}(b)$,
$a,b \in U$, where we are using the bar $\bar{~~}$ for images under the
homomorphism from $E$ to $\bar{E}$.
As in \cite{Asch}, Theorem~23.10 we define a
quadratic form $Q: \bar{E} \ra \FF_2$ by
$$
Q( \bar{g} ) = \left\{
\begin{array}{ll}
0 & \mbox{if $g^2 = + I$} \\
1 & \mbox{if $g^2 = - I$}
\end{array}
\right.
$$
for $\bar{g}\in \bar{E}$, where $g \in E$ is any preimage of $\bar{g}$,
and so $Q( \bar{X} (a) \bar{Y} (b)) = a \cdot b$.

The associated alternating bilinear form $B: \bar{E} \times \bar{E} \ra \FF_2$
is given by
$$
B( \bar{g}_1, \bar{g}_2 ) = Q( \bar{g}_1 +
\bar{g}_2 ) + Q( \bar{g}_1 ) + Q( \bar{g}_2 )~,
$$
for $\bar{g}_1 , \bar{g}_2 \in \bar{E}$, and so
\beql{eq5}
B( \bar{X} (a) \bar{Y} (b)~,~ \bar{X} (a')
\bar{Y} (b')) = a \cdot b' + a' \cdot b~.
\eeq
Then $(\bar{E},Q)$ is an orthogonal vector space of type $\Om^+ (2i,2 )$
and maximal Witt index (cf. \cite{Dieu}).

The normalizer of $E$ in $\sO$ is a
certain Clifford group $L$, of order
$$
2^{i^2 + i +2}~~ (2^i -1)~~ \Pi_{j=1}^{i-1} (4^j -1 )
$$
(cf. [CCKS], Section~2).
$L$ is generated by $E$, all permutation matrices
$G(A,a) \in \sO : e_u \ra e_{Au+a} $,
$u \in U$, where
$A$ is an invertible $i \times i$ matrix over $\FF_2$ and $a \in U$, and the further matrix
$H = (H_{u,v} ) $,
$H_{u,v} = 2^{- i/2} (-1)^{u \cdot v}$,
$u,v \in U$.

The group $L$ acts on $E$ by conjugation, fixing the center, and
so also acts on $\bar{E}$.
In fact $L$ acts on $\bar{E}$ as the orthogonal
group $O^+ (2i,2)$
(\cite{CCKS}, Lemma~2.14).

This Clifford group $L$ has arisen in several different contexts,
providing a link between the present problem,
the Barnes-Wall lattices (see \cite{BRWI},
\cite{BRWII}, \cite{grass2}, \cite{Wall}), the construction of
orthogonal spreads and Kerdock sets \cite{CCKS}, and the
construction of quantum error-correcting codes \cite{BDSW},
\cite{QC2}.
It also occurs in several purely group-theoretic contexts -- see
\cite{CCKS} for references.

The connection with quantum computing arises because if certain
conditions are satisfied the invariant
subspaces mentioned in Theorem~1 form good
quantum-error-correcting codes \cite{QC1}, \cite{QC2}.
\section{The construction from totally singular subspaces}
\label{sec3}
\hsp
A subspace $\bar{S} \subseteq \bar{E}$ is {\em totally singular}
if $Q( \bar{g}) = 0$ for all
$\bar{g} \in \bar{S}$.
Then $\dim \bar{S} \leq i$,
and if $\dim \bar{S} = i$ then $\bar{S}$ is
{\em maximally totally singular}.
It follows from (\ref{eq5}) that the preimage
$T \subseteq E$ of a maximally totally singular space $\bar{T}$ is an
abelian subgroup of $E$, of order $2^{i+1}$.
$T$ contains $-I$, and has $2^{i+1}$ linear characters,
associated with $2^i$ mutually perpendicular
1-dimensional invariant subspaces forming a
coordinate frame $\sF (T) \subset V$
(\cite{CCKS}, Lemma~3.3).

Since $L$ acts as $O^+ (2i,2)$ on $\bar{E}$,
$L$ takes
any ordered pair of maximally
totally singular subspaces that meet in $\{0\}$ to $X$ and $Y$ respectively.
The corresponding coordinate frames in $V$ are
\beql{eq6}
\sF (X) = \{ e_v^\ast = \df{1}{2^{i/2}} \sum_{u \in U}
(-1)^{u \cdot v} e_u : v \in U \}
\eeq
and
\beql{eq7}
\sF (Y) = \{e_u : u \in U \}~,
\eeq
respectively.

If $\bar{S} \subseteq \bar{T}$ has dimension $k$, its
preimage $S \subseteq E$ has $2^{k+1}$ linear characters, and
$2^k$ distinct invariant subspaces, each of which is spanned by
$2^{i-k}$ of the vectors in $\sF (T)$.

We can now state our first main construction for Grassmannian packings.
\paragraph{\bf Theorem 1.}
{\em Given $k$, with $0 \leq k \leq i-1$, the set of all invariant subspaces of the preimages $S$
of all $(i-k)$-dimensional totally singular subspaces $\bar{S}$ of $\bar{E}$
is a packing of $N$ planes in $G(2^i , 2^k )$ with minimal distance
$d = 2^{(k-1)/2}$, where}
$$
N = 2^{i-k} \left[
\begin{array}{l}
i \\
k
\end{array}
\right] \prod_{j=k}^{i-1} (2^j +1) ~,
$$
{\em and}
$$
\left[
\begin{array}{l}
i \\
k
\end{array}
\right] = \df{(2^i -1) \ldots (2^{i-k+1} -1)}{(2^k -1) \ldots (2-1 )}
$$
{\em is a Gaussian binomial coefficient.}
\paragraph{\bf Proof.}
There are
$$
\left[
\begin{array}{l}
i \\
k
\end{array}
\right] \prod_{j=k}^{i-1} (2^j +1)
$$
choices for $\bar{S}$ (\cite{BCN}, Lemma~9.4.1)
and each $\bar{S}$
yields $2^{i-k}$ planes.

We compute the distance between two planes from (\ref{eq2}).
Suppose $P_j$ is a $2^k$-dimensional invariant subspace of $V$
corresponding to the character $\chi_j$ of the subgroup
$S_j \subseteq E$, for $j = 1,2$.
We may assume $-I \in S_1 \cap S_2$
and $\chi_1 (-I) = \chi_2 (-I) = -1$.
Then
$$
\Pi_j = \df{1}{|S_j |} \sum_{g \in S_j} \chi_j (g) g \in \sO
$$
is the orthogonal projection onto $P_j$.
Also
\begin{eqnarray}
\mbox{trace}~ (\Pi_1 ~ \Pi_2 ) & = &
\df{1}{|S_1 |~|S_2|}
\sum_{g_1 \in S_1} ~
\sum_{g_2 \in S_2 }
\chi_1 (g_1) \chi_2 (g_2 ) ~\mbox{trace}~
(g_1 g_2 ) \nonumber \\
~~~ \nonumber \\
  & = & \df{1}{|S_1 |~|S_2 |} \sum_{g_1 \in S_1 \cap S_2} ~
\sum_{g_2 = \pm g_1^{-1}}
\chi_1 (g_1 ) \chi_2 (g_2 ) ~\mbox{trace}~ (g_1 g_2 ) \nonumber \\
~~~ \nonumber \\
  & = & \df{2}{|S_1 |~|S_2 |} \sum_{g_1 \in S_1 \cap S_2}
\chi_1 (g_1 ) \chi_2 (g_1^{-1} ) ~\mbox{trace}~ (I) \nonumber \\
~~~ \\
 & = & \left\{
\begin{array}{ll}
\df{2|S_1 \cap S_2 | 2^i}{|S_1 |~|S_2 |} ~, & \mbox{if $\chi_1 = \chi_2$
on $S_1 \cap S_2$} \nonumber \\
~~~ \nonumber \\
0 ~, & \mbox{otherwise}~.
\end{array}
\right.
\end{eqnarray}
This implies from (\ref{eq2}) that
\beql{eq9}
d^2 (P_1 ,~P_2 ) = \left\{
\begin{array}{ll}
2^k - \df{| \overline{S_1 \cap S_2 |}}{|S_1 |~|S_2|} 2^i & \mbox{if $\chi_1 = \chi_2$ on $S_1 \cap S_2$} \nonumber \\
~~~ \nonumber \\
2^k & \mbox{otherwise}~.
\end{array}
\right.
\eeq
So if $S_1 = S_2$ and $\chi_1 \neq \chi_2$ the planes are
orthogonal and at distance
$2^{k/2}$; otherwise $S_1 \neq S_2$,
$| \overline{S_1 \cap S_2} | \leq 2^{i-k-1}$, and
their distance satisfies
$$
d^2 \geq 2^k - 2^{k-1} = 2^{k-1} ~,
$$
as claimed.
\hfill $\btd$

The principal angles $\tht_1 , \ldots , \tht_{2^k}$
between any two planes in the packing may be found by a
similar calculation, using the fact that the
singular values of $\Pi_1 \Pi_2$ are $\cos^2 \tht_1 , \ldots , \cos^2 \tht_{2^k}$
together with $2^i - 2^k$ zeros.
It turns out that the principal angles are either all equal to $\pi/2$, or
else $N_1$ of them are equal to
arccos $2^{-r/4}$ and
$2^k -N_1$ are equal to $\pi/2$, where $r$ is the rank of $Q$
on $\bar{S}_1 \cup \bar{S}_2$ and
$$
N_1 = 2^{2k-i+r/2} |
\bar{S}_1 \cap \bar{S}_2 |~.
$$
\paragraph{\bf Examples.}
Taking $k=0$ in the theorem we obtain a packing of
$$
(2+2)(2^2 +2) \ldots (2^i +2)
$$
lines in $G(2^i , 1)$ with minimal angle $\pi/4$
(as in \cite{grass2}).
These are the lines defined by the minimal vectors in the
$2^i$-dimensional Barnes-Wall lattice together with their images
under $H$
(cf. \cite{SPLAG}, p.~151).$^4$\footnotetext[4]
{The group of the Barnes-Wall lattice in dimension $2^i$ is a subgroup
$G$ of index~2 in $L$.
This lattice can therefore be constructed by taking unit
vectors along the coordinate
frame $\sF (Y)$, forming their images under $G$ and then
taking their integral closure.}

With $i = 2,~k=1$ and $i = 3,~k=2$
we obtain two important special cases:
18~planes in $G(4,2)$ and 70~planes in $G(8,4)$
(cf. \cite{grass1}, \cite{grass2}).
More generally, when $k = i -1$
we obtain the packing of
$$
f(i) = 2(2^i -1)(2^{i-1} +1 )
$$
planes in $G(2^i , 2^{i-1})$ with $d^2 = 2^{i-2}$ that
is the main result of \cite{grass2}.
These packings meet the orthoplex bound of (\ref{eq4}) and are
therefore optimal.
(We do not know if any of the other examples are optimal.
Even if not optimal as Grassmannian packings, they may be
optimal subject to constraints on the spectrum of
distances between planes -- compare Prop.~3.12 of [CKSS].)
An explicit recursive construction for the special case
$k = i-1$ is given in \cite{grass2}.

For $k =1$ and $k = i-2$ we obtain two further sequences of
packings whose existence was conjectured in \cite{grass2}.

The construction given in the theorem can be restated in an equivalent but more
explicit way as follows.
Let $P_0$ be the $2^k$-dimensional plane spanned by the
coordinate vectors $e_u$, where $u \in U$ is of the form
$00 \ldots 0 \ast \ldots \ast$,
with $i-k$ initial zeros.
Then the packing consists of all the images of $P_0$ under the group $L$.

The parameters of all the packings obtained from the
theorem in dimensions up to 128 can be seen in Table~I below.

Many other packings can be obtained from the images under $L$ of other
starting planes, and still more by replacing $L$ by smaller groups.
We mention just one such family.
With the $m=2^i$ coordinate vectors $e_u \in V$ arranged in
lexicographic order, let $\sL (m,n)$ denote the packing
in $G(m,n)$ obtained from the images under $L$ of the plane spanned
by the first $n$ coordinate vectors.
The $\sL (2^i , 2^k )$ are the packings described in Theorem~1.
In particular, $\sL (2^i , 2^{i-2} )$
contains $\slashedfrac{1}{3}~ f(i) f(i-1)$ planes and has $d^2 = 2^{i-3}$.
The numerical evidence (see the entries
marked ``(1a)'' in Table~I)
suggests that $\sL (2^i , 3.2^{i-3} )$
contains
$\slashedfrac{1}{3}~ f(i) f(i-1) f(i-2)$ planes and
has $d^2 = 2^{i-4}$,
and that
$\sL (2^i , 5.2^{i-4} )$
contains
$\slashedfrac{1}{3}~ f(i) f(i-1) f(i-2) f(i-3)$ planes and has
$d^2 = 2^{i-5}$.
\section{Spreads and clique-finding}
\label{sec4}
\hsp
The packings constructed in Section~\ref{sec3} contains very
large numbers of planes.
Smaller packings can be obtained by using only some of the totally
singular spaces.
\paragraph{\bf Theorem 2.}
{\em Suppose a set of $M$ totally singular $(i-k)$-subspaces $\bar{S}$ of $\bar{E}$ can be
found such that any pair intersect in a space of dimension at most $\ell$.
Then the set of invariant subspaces of all the preimages $S \subseteq E$ is a packing of $2^{i-k} M$ planes in $G(2^i ,~ 2^k )$ with minimal distance
satisfying }
\beql{eq10}
d^2 \geq 2^k - 2^{2k+\ell -i}~.
\eeq
\paragraph{\bf Proof.}
The bound on the minimal distance follows from (\ref{eq9}).
Equality holds in (\ref{eq10}) if and only if some pair of the
subspaces intersect in a subspace of dimension exactly $\ell$.
\hfill $\btd$
\paragraph{\bf Examples.}
(a)
An {\em orthogonal spread} \cite{CCKS}, \cite{STK}, \cite{OST} is a partition of the nonzero
totally singular points of $\bar{E}$ into $2^{i-1} +1$ totally singular $i$-spaces.
Such a partition exists if and only if $i$ is even (the construction is
closely related to Kerdock codes), and then Theorem~2 applies
with $M = 2^{i-1} +1$,
$k = \ell = 0$,
producing $2^i (2^{i-1} +1)$ lines in $G(2^i ,~1)$ with minimal angle
arccos $2^{-i/2}$.

(b) More generally, a {\em spread}
(\cite{Hirsch}, Theorem~4.1.1) in a projective space $PG(s,q)$ is a partition of the
points into copies of $PG(r,q)$, and exists if and only if $r+1$ divides $s+1$.
Suppose now $i$ is even and $j$ divides $i$.
If we take every totally singular $i$-space in an orthogonal spread and partition
the nonzero points into copies
of a $PG(j-1,2)$,
using a spread, we obtain $M = (2^i -1)/(2^j -1)$ totally
singular $j$-spaces in $\bar{E}$ which meet only in the zero vector.
This produces a packings of
$2^j (2^{i-1} +1)(2^i -1)/(2^j -1)$ planes in $G(2^i , 2^{i-j} )$
with $d^2 = 2^{i-j} - 2^{i-2j}$.

In particular, because $i$ is even we can always take $j = 2$,
obtaining $4(2^{i-1} +1)(2^i -1)/3$ planes in
$G(2^i ,~ 2^{i-2} )$ with
$d^2 = 3.2^{i-4}$.
These packings meet the orthoplex bound of (\ref{eq4}).

(c) When the general constructions in (a) and (b) are not applicable,
or do not give the desired parameters, we may always resort to
clique-finding.
We form a graph whose nodes represent all totally singular $(i-k)$-spaces
$\bar{S} \subseteq \bar{E}$, with an
edge joining two nodes if the corresponding spaces
intersect in a space of dimension at most $\ell$, and search for a
maximal clique.
Theorem~2 gives the parameters of the resulting packing.

For example, when $i = 3, k = 1, \ell = 0$, the graph on
2-spaces has 105~nodes and contains maximal cliques of size~10,
producing packings of 40~planes in
$G(8,2)$
with $d^2 = 1.5$.
These are suboptimal however, since packings of 44~planes in
$G(8,2)$ with $d^2 = 1.5$ were obtained in \cite{grass1}.

For $i = 4$, $k = 1 , \ell = 0$, the graph on 3-spaces has 2025~nodes and
contains cliques of size~17
(which is probably maximal),
leading to packings of 136~planes in $G(16,2)$ with $d^2 = 1.75$.
For $\ell =1$ the graph contains cliques of size at
least 130, giving 1040~planes in $G(16,2)$ with $d^2 = 1.5$.
We do not know if these are optimal.

Instead of orthogonal spreads in real space, we can also make use of
their analogues in complex or quaternionic space.
Since the packings obtained do not seem especially good we give only a summary.

(d)
A {\em symplectic spread} is the complex analogue of an
orthogonal spread, and leads to a family of $2^{i-1} (2^{i-1} +1)$
vectors in $\CC^{2^{i-1}}$
whose angles are $\pi /2$ or arccos
$2^{-(i-1)/2}$,
for $i \geq 2$ (\cite{CCKS}, Theorem~5.6).
This produces packings of $2^{i-1} (2^{i-1} +1)$
planes in $G(2^i , 2)$ with $d^2 = 2 -2^{-(i-2)}$.

(e)
In a similar way, Kantor \cite{QLS} constructs a family of $2^{i-2} (2^{i-1} +1)$
lines in quaternionic space of dimension $2^{i-2}$ whose
angles are $\pi/2$ or arccos $2^{-(i-2)/2}$, for all
even $i \geq 4$.
This produces packings of $2^{i-2} (2^{i-1} +1)$ planes in $G(2^i , 4)$ with $d^2 = 4-2^{-(i-2)} )$.
\section{An infinite family of packings meeting the simplex bound}
\label{sec5}
\hsp
A packing of 28 planes in $G(7,3)$ meeting the simplex bound of
(\ref{eq3}) was given in \cite{grass1}.
This may be generalized as follows.

Let $p$ be a prime which is either 3 or congruent to $-1$ modulo~8, so that a Hadamard matrix $H$ of order $(p+1)/2$ exists.
Let $Q = \{q_1 , \ldots , q_{(p-1)/2} \}$ denote the
nonzero quadratic residues modulo $p$, and
$R = \{r_1 , \ldots , r_{(p-1)/2} \}$ the nonresidues.
The entries of $H$ will be denoted by $H_{s,t}$,
for $0 \leq s,t \leq (p-1)/2$, where we assume $H_{s,0} = H_{0,t} = 1$ for
all $s,t$.
We will construct a packing in $G(p, \frac{p-1}{2} )$.
Let $e_s (0 \leq s \leq p-1)$ be coordinate vectors
in $\RR^p$,
let $C = (1+ \sqrt{p+2} )/\sqrt{p+1}$, and fix an element $k \in R$.
Let $P_t (0 \leq t \leq (p-1)/2)$ be the
$(p-1)/2$-dimensional plane spanned by the vectors
$$
e_{q_s} ~~+~~ H_{st}~~C~~e_{kq_s}~, \quad 1 \leq s \leq \df{p-1}{2} ~.
$$
For each $P_t$ we obtain $p-1$ further planes by applying the cyclic permutation
of coordinates $e_i ~\ra~ e_{i+1 (\bmod~p )}$, for
a total of $p(p+1)/2$ planes.

\paragraph{\bf Theorem 3.}
{\em The above construction produces a packing of $p(p+1)/2$ planes in
$G(p, \frac{p-1}{2} )$ in which the distance between every pair of
planes satisfies}
$$
d^2 = \df{(p+1)^2}{4(p+2)} ~.
$$
{\em This packing meets the simplex bound of (\ref{eq3}) and is therefore optimal.}
\paragraph{\bf Sketch of proof.}
The principal angles between two planes in the same orbit under the
cyclic shift, for example $P_0$ and $P_1$, are
$0~~ ( \frac{p-3}{4}$ times) and
arcsin $2C/(1+C^2 )~~$
$( \frac{p+1}{4}$ times), and so
\beql{eq50}
d^2 (P_0 , P_1 ) = \df{p+1}{4}
\df{4C^2}{(1+C^2 )^2} =
\df{(p+1)^2}{4(p+2)} ~.
\eeq

If the two planes are in different orbits, it is best to compute the
corresponding projection matrices and then use (\ref{eq2}) to compute the
distance, making use of Perron's theorem \cite{Perron} on
quadratic residues.
Again the distance is given by (\ref{eq50}).
We omit the details of this calculation.
\hfill $\btd$

For example, when $p = 7$ and $C = \sqrt 2$,
taking $k = 3$ we find that the planes $P_0 , \ldots , P_3$
are generated by
$$
\Blb{.3}
\begin{array}{c c c c c c c}
 0 & 1 & 2 & 3 & 4 & 5 & 6 \\[.05in]
0 & 1 & 0 & \pm C & 0 & 0 & 0 \\
0 & 0 & 1 & 0 & 0 & 0 & \pm C \\
0 & 0 & 0 & 0 & 1 & \pm C & 0 \\
~ & ~ & ~ & ~ & ~ & ~ & ~
\end{array}
\Brb{.3}~,
$$
where the product of the signs is $+1$ (as given in \cite{grass1}).
The full set of 28~planes is obtained by
cycling the seven coordinates.
Changing $k$ to 6 we obtain a packing with the same distances but in
which some of the principal angles have changed, showing that the
packings of Theorem~3 are not unique.

Packings obtained from Theorem~3 are labeled ``(3)''
in Table~I.
\section{A table}
\label{sec6}
\hsp
Table~I lists the parameters of the packings constructed in this paper
in dimensions up to 128.
When better packings with same parameters were given in
\cite{grass1},
these are also mentioned.
In the last column, ``(1)'' refers to Theorem~1,
``(1a)'' to the packings described at the end of Section~3,
``(2a)''$, \ldots ,$
``(2e)'' to the examples following Theorem~2, and ``(3)'' to
Theorem~3.

We must stress however that a very large number of
other packings are known, especially in dimensions up to 16:
see the constructions and tables in
\cite{grass1} and \cite{grass0}.

\clearpage
\begin{center}
{\em Table~I.}~~Parameters of Grassmannian Packings constructed in this paper 

\vspace{.5\baselineskip}

\begin{tabular}{lrrrl}
$m$ & $n$ & $N$ & $d^2$ & source   \\ \hline
3 & 1 & 6 & 4/5 & (3) \\
 & & & & \\
4 &    1 &    12&    3/4 &   (2a) \\
4 &    1 &    24 &     1/2 &    (1) \\
4 &    1 &   24 &   0.5182$\ldots$ &  \cite{grass1} \\
4 &    2 &    6 &    1 &    (2d) \\
4 &   2 &   6 &    6/5 &  \cite{grass1} \\
4 &    2 &    18 &   1   & (1) \\
 & & & & \\
7 & 3 & 28 & 16/9 & (3) \\
 & & & & \\
8 &    1 &    240 &    1/2 &    (1) \\
8 &    2 &    20 &    3/2 &    (2d) \\
8 &    2  & 20 &  1.5789$\ldots$ &  \cite{grass1} \\
8 &    2 &    40 &    3/2 &   (2c) \\
8 &    2 & 44 &  3/2 &  \cite{grass1} \\
8 &    2 &    420 &    1 &    (1) \\
8 & 3 & 1680 & 1/2 & (1a) \\
8 &    4 &    70 & 2 &   (1) \\
 & & & & \\
16 &    1 &  144 &    15/16 &   (2a) \\
16  &   1 &  4320 &    1/2 &    (1) \\
16 &    2 &    72 &    7/4 &    (2d) \\
16 &    2 &   136 &   7/4 &   (2c) \\
16 &    2 &   1040 &   3/2 &   (2c) \\
16 &    2 &    16200 &   1 &   (1) \\
16 & 3 & 151200 & 1/2 & (1a) \\
16 &  4  &   72 &    15/4  &   (2e) \\
16 &  4  &   180 &   3 &    (2b) \\
16 &  4  &   6300 &   2 &   (1) \\
16 & 5 & 453600 & 1/2 & (1a) \\
16 & 6 & 113400 & 1 & (1a) \\
16 & 7 & 64800 & 1/2 & (1a) \\
16 &  8   & 270 &   4  &   (1) \\
 & & & & \\
23 & 11 & 276 & 144/25 & (3) \\
 & & & & \\
31 & 15 & 496 & 256/33 & (3) \\
 & & & & \\
32 &  1 &    146880 &   1/2 &    (1) \\
32  &  2 & 272 &    15/8 &    (2d) \\
32 &  2 &  1138320 &   1 &   (1) \\
32 &  4  &   948600 &    2 &  (1) \\
32 &  8 &    94860 &    4  &   (1) \\
32 &  16 &  1054  &   8 &    (1) \\
\end{tabular}
\end{center}
\clearpage
\begin{center}
{\em Table~I. cont'd.}~~Parameters of Grassmannian Packings constructed in this paper

\vspace{1.5\baselineskip}

\begin{tabular}{lrrrl}
$m$ & $n$ & $N$ & $d^2$ & source \\ \hline
47 & 23 & 1128 & 576/49 & (3) \\
 & & & & \\
64  &    1 &    2112 &    63/64  &   (2a) \\
64  &    1 &    9694080 &    1/2 &    (1) \\
64  &    2 &    1056 &    31/16 &  (2d) \\
64  &    2 &    152681760 &    1 &   (1) \\
64  &    4  &   1056 &    63/16 &  (2e) \\
64  &    4  &   262951920 &    2 &  (1) \\
64  &    8 &    2376 &   7  & (2b) \\
64 & 8 & 56346840 & 4 & (1) \\
64  &    16 &  2772 &    12 &   (2b) \\
64  &    16 &  1460844  &   8 &  (1) \\
64  &    32 &  4158 &    16 &  (1) \\
 & & & & \\
71 & 35 & 2556 & 1296/73 & (3) \\
 & & & & \\
79 & 39 & 3160 & 1600/81 & (3) \\
 & & & & \\
103 & 51 & 5356 & 2704/105 & (3) \\
 & & & & \\
127 & 63 & 8128 & 4096/129 & (3) \\
 & & & & \\
128 &   1 &   1260230400 &   1/2 &    (1) \\
128 &   2 &   4160 &    63/32 &  (2d) \\
128 &   2 &   40012315200 &    1  & (1) \\
128 &   4  &   140043103200 &   2 &   (1) \\
128 &   8   & 62019088560 &    4  &   (1) \\
128 &   16 &  3445504920 &    8 &  (1) \\
128 &   32 &  22882860 &    16 &  (1) \\
128 &   64  &   16510 &   32 &  (1) \\
\end{tabular}
\end{center}

\end{document}